\documentclass[12pt,leqno]{amsart}

\usepackage[dvips]{graphicx}
\usepackage{float}
\usepackage{amssymb,amsmath,amsthm,amscd}
\usepackage{graphics}
\usepackage{mathrsfs}
\usepackage{enumerate}
\usepackage{color}
\usepackage{verbatim}

\usepackage[bookmarks=true,bookmarksnumbered=true,bookmarkstype=toc=true,pdftitle={\@title},pdfauthor={\@author},pdfstartview={FitBH -32768}]{hyperref}

\usepackage[all]{xy}
\numberwithin{equation}{section}
\newcommand{\nc}{\newcommand}
\usepackage{graphicx}
\usepackage{tikz}

\def\node#1#2{\overset{#1}{\underset{#2}{\circ}}}

\def\Llleftarrow{%
\lower2pt\hbox{\begingroup
\tikz
\draw[shorten >=0pt,shorten <=0pt] (0,3pt) -- ++(-1em,0) (0,1pt) -- ++(-1em-1pt,0) (0,-1pt) -- ++(-1em-1pt,0) (0,-3pt) -- ++(-1em,0) (-1em+1pt,5pt) to[out=-105,in=45] (-1em-2pt,0) to[out=-45,in=105] (-1em+1pt,-5pt);
\endgroup}
}

\DeclareMathOperator{\WORD}{\mathsf{word}}
\DeclareMathOperator{\SHAPE}{\mathsf{sh}}

\DeclareMathOperator{\CONT}{\mathsf{cont}}
\DeclareMathOperator{\RSK}{\mathsf{RSK}}
\DeclareMathOperator{\MUMAP}{\mathsf{M}}
\DeclareMathOperator{\WT}{wt}
\DeclareMathOperator{\SWAP}{\mathsf{swap}}
\DeclareMathOperator{\SHIFT}{\mathsf{shift}}
\DeclareMathOperator{\BASE}{\mathsf{base}}

\DeclareMathOperator{\ROOF}{\mathsf{roof}}
\DeclareMathOperator{\ST}{\mathsf{ST}}
\DeclareMathOperator{\RST}{\mathsf{RST}}
\DeclareMathOperator{\SST}{\mathsf{SST}}
\DeclareMathOperator{\RSST}{\mathsf{CSPP}}
\DeclareMathOperator{\PP}{\mathsf{PP}}

\newcommand{\JOKEN}{\ast}
\newcommand{\EUPHI}{\phi}

\newcommand{\GEE}{\mathfrak{g}}
\newcommand{\GEEO}{\check{\mathfrak{g}}}
\newcommand{\AO}{\check{A}}
\newcommand{\IO}{\check{I}}
\newcommand{\DEO}{\check{\delta}}
\newcommand{\LAMBDAO}{\check{\Lambda}_0}
\newcommand{\LAMBDA}{\check{\Lambda}}
\newcommand{\ALPHAO}[1]{\check{\alpha}_{#1}}
\newcommand{\WTO}{\check{\WT}}

\newcommand{\HA}{\mathfrak{h}}

\newcommand{\diagautom}{\omega}
\newcommand{\SIGN}[1]{\mathsf{sign}_{#1}}

\DeclareMathOperator{\DDD}{\mathsf{D}}

\newcommand{\HSL}[1]{\widehat{\mathfrak{sl}_{#1}}}

\newcommand{\LAM}[3]{\nu^{#1}_{#2,#3}}
\newcommand{\MYMULT}[3]{m_{#3}(#1,#2)}
\newcommand{\MYWEIGHT}[2]{P_{#2}(#1)}
\newcommand{\SPECHT}[2]{D^{#2}_{\mathbb{F}_{#1}}}
\newcommand{\MYMAX}[2]{\max_{#2}(#1)}

\newcommand{\TOKUMEI}{T}
\newcommand{\ANNO}{\Phi}
\newcommand{\ANNOO}{\Psi}

\newcommand{\emptypartition}{\varnothing}

\newcommand{\IRR}{\mathsf{Irr}}

\newcommand{\PAR}{\mathsf{Par}}

\newcommand{\RPAR}[1]{\mathsf{RPar}_{#1}}
\newcommand{\RPARR}[2]{\mathsf{RPar}^{#2}_{#1}}
\newcommand{\CRPAR}[1]{\mathsf{Par}^{\textrm{$#1$-core}}}

\newcommand{\C}{{\mathbb C}}

\newcommand{\Q}{\mathbb {Q}}

\newcommand{\Z}{{\mathbb Z}}

\newcommand{\MP}[1]{\mathcal{P}_{#1}^{+}}

\newcommand{\corr}{\mathbb{F}}

\newcommand{\TRANS}[1]{{}^{\textrm{tr}}{#1}}

\theoremstyle{plain}
\newtheorem{lemma}{Lemma}[section]
\newtheorem{prop}[lemma]{Proposition}
\newtheorem{theorem}[lemma]{Theorem}
\newcommand{\Prop}{\begin{prop}}
\newcommand{\enprop}{\end{prop}}
\newcommand{\Lemma}{\begin{lemma}}
\newcommand{\enlemma}{\end{lemma}}
\newcommand{\Th}{\begin{theorem}}
\newcommand{\enth}{\end{theorem}}
\newtheorem{corollary}[lemma]{Corollary}
\newcommand{\Cor}{\begin{corollary}}
\newcommand{\encor}{\end{corollary}}
\newtheorem{definition}[lemma]{Definition}
\newtheorem{question}[lemma]{Question}
\newcommand{\Question}{\begin{question}}
\newcommand{\enquestion}{\end{question}}
\newcommand{\Def}{\begin{definition}}
\newcommand{\edf}{\end{definition}}

\theoremstyle{definition}
\newtheorem{remark}[lemma]{Remark}

\newtheorem{conjecture}[lemma]{Conjecture}
\nc{\Rem}{\begin{remark}}
\nc{\enrem}{\end{remark}}

\newcommand{\Con}{\begin{conjecture}}
\newcommand{\encon}{\end{conjecture}}

\nc{\emprule}[1]{\rule{#1}{0pt}}
\newcommand{\g}{{\mathfrak{g}}}

\newcommand{\isoto}[1][]{\mathop{\xrightarrow[#1]{\rule{0pt}{.9ex}{\raisebox{-.6ex}[0ex][-.7ex]{$\mspace{4mu}\sim\mspace{3mu}$}}}}}

\renewcommand{\hom}{\operatorname{\it \mathscr{H}\kern-.25em om}}

\newcommand{\M}{{\mathscr M}}
\newcommand{\N}{{\mathbb{N}}}
\newcommand{\eq}{\begin{eqnarray}}
\newcommand{\eneq}{\end{eqnarray}}

\newcommand{\eqn}{\begin{eqnarray*}}
\newcommand{\eneqn}{\end{eqnarray*}}

\newenvironment{tenumerate}{
  \begin{enumerate}
  
  }{\end{enumerate}}

\nc{\bnum}{\begin{enumerate}[{\rm (i)}]}
\nc{\enum}{\ee}
\nc{\bna}{\begin{enumerate}[{\rm (a)}]}
\nc{\bnA}{\begin{enumerate}[{\rm (A)}]}
\nc{\bnX}{\begin{enumerate}[{\rm (X)}]}
\nc{\ena}{\ee}

\newcommand{\on}{\operatorname}
\newcommand{\bni}{\begin{tenumerate}}
\newcommand{\eni}{\end{tenumerate}}

\newcommand{\QED}{\end{proof}}
\newcommand{\Proof}{\begin{proof}}

\newcommand{\ba}{\begin{array}}
\newcommand{\ea}{\end{array}}

\newcommand{\hs}{\hspace*}

\newcommand{\eqsub}{\begin{subequations}\begin{eqnarray}}
\newcommand{\eneqsub}{\end{eqnarray}\end{subequations}}

\newcommand{\ol}{\overline}

\newcommand{\A}{\mathscr{A}}

\nc{\la}{\lambda}
\nc{\lam}{\lambda}
\nc{\U}[1][\g]{U_q(#1)}
\nc{\te}{\tilde{e}}
\nc{\tei}{\tilde{e}_i}
\nc{\tf}{\tilde{f}}
\nc{\KF}[1]{{\tilde{f}}_{#1}}
\nc{\tU}{\widetilde U_q(\g)}
\nc{\tE}{\tilde{E}}
\nc{\tF}{\tilde{F}}

\nc{\BZ}{{\mathbb{Z}}}
\nc{\al}{\alpha}
\nc{\qs}{{q}}
\nc{\lan}{\langle}
\nc{\ran}{\rangle}
\nc{\re}{{\mathrm{re}}}
\nc{\wt}{\operatorname{wt}}
\nc{\Uf}[1][\g]{U^-_q(#1)}
\nc{\Ue}{U^+_q(\g)}
\nc{\eps}{\varepsilon}
\nc{\vphi}{\varphi}
\nc{\sphi}{\varphi^*}
\nc{\seps}{\varepsilon^*}

\nc{\nn}{\nonumber}
\def\max{{\mathop{\mathrm{max}}}}
\nc{\vp}{\varpi}
\nc{\cls}{{\operatorname{cl}}}

\nc{\Wt}{{\operatorname{Wt}}}
\nc{\Us}{U'_q(\g)}
\nc{\La}{\Lambda}
\nc{\ro}{{\rm(}}
\nc{\rf}{{\rm)}}
\nc{\norm}{{\mathrm{norm}}}
\nc{\qbox}{\quad\mbox}
\nc{\braid}{{\mathfrak{B}}}
\nc{\Ad}{\operatorname{Ad}}
\nc{\Aut}{\operatorname{Aut}}
\nc{\dt}[1]{\tilde{\tilde #1}}
\nc{\Sn}{S^{{\mathrm{norm}}}}
\nc{\aff}{{\mathrm{aff}}}
\nc{\rk}{{\mathrm{rk}}}
\nc{\tQ}{\widetilde{Q}}
\nc{\tP}{\widetilde{P}}
\nc{\tW}{\widetilde{\mathscr{W}}}
\nc{\Dyn}{\mathrm{Dyn}}
\nc{\tD}{\widetilde{\Delta}}
\nc{\height}{{\operatorname{ht}}}
\nc{\bl}{\bigl}
\nc{\br}{\bigr}
\nc{\Hecke}{\mathrm{H}}
\nc{\HB}{\Hecke^{\mathrm{B}}}
\nc{\K}{\mathrm{K}}
\newcommand{\scbul}{{\,\raise1pt\hbox{$\scriptscriptstyle\bullet$}\,}}
\nc{\vac}{{\phi}}
\nc{\be}{\begin{enumerate}}
\nc{\ee}{\end{enumerate}}
\nc{\low}{{\mathrm{low}}}
\nc{\upper}{{\mathrm{up}}}
\nc{\Zodd}{\Z_{\mathrm{odd}}}
\nc{\Ft}[1][n]{\mathbb{P}\mathrm{ol}_{#1}}
\nc{\Ftf}[1][n]{\widetilde{\mathbb{P}\mathrm{ol}}_{#1}}
\nc{\KB}{\on{K}^{\mathrm{B}}}
\nc{\Res}{\on{Res}}
\nc{\Fc}[1][{n,m}]{\mathbf{F}_{#1}}
\nc{\tphi}{\widetilde{\varphi}}
\nc{\CO}{\mathscr{O}}
\nc{\CK}{\mathscr{K}}
\nc{\disc}{\mathfrak{d}}
\nc{\tr}{\on{tr}}
\nc{\Gb}{\mathfrak{b}}
\nc{\Gh}{\mathfrak{h}}
\nc{\ga}{\mathfrak{a}}
\nc{\stable}{\mathrm{stable}}
\nc{\X}{\mathfrak{X}}
\nc{\Hilb}{\mathrm{Hilb}}
\nc{\W}{\ensuremath{\mathscr{W}}}
\nc{\Ws}{\ensuremath{\rm W}}
\nc{\opp}{{\on{opp}}}

\nc{\corps}{{\mathbf{k}}}
\nc{\cor}{{\mathbf{k}}}
\nc{\h}{\mathrm{\hslash}}
\nc{\fL}[1][{\h}]{\C(\mspace{-1mu}(#1)\mspace{-1mu})}
\nc{\ad}{\mathrm{ad}}
\newcommand{\Endm}{\operatorname{\mathscr{E}\kern-.1pc\mathit{nd}}}
\newcommand{\Endomo}{\operatorname{\mathscr{E}\kern-.1pc\mathit{nd}}}
\nc{\bc}{\bar{\corps}}
\nc{\reg}{{\mathrm{reg}}}
\nc{\ysq}{\mathbf{y}^2}
\nc{\CH}{\mathsf{char}}

\nc{\sketch}{\Proof}
\nc{\Gm}{\mathbb{G}_{\mathrm{m}}}
\nc{\hGm}{\hat{\mathbb{G}}_{\mathrm{m}}}
\nc{\ug}{\widehat{\mathrm{G}}_{\mathrm{m}}}

\nc{\tL}{\widetilde{\mathscr{L}}}
\nc{\Fr}{\mathcal{F}}
\nc{\E}{\mathcal{E}}

\nc{\ord}{\on{ord}}
\nc{\bM}{\overset{\hs{1.5ex}\rule[-.08ex]{1.8ex}{.08ex}}{\M}}
\nc{\romano}{\mathrm{o}}
\nc{\into}{\hookrightarrow}
\nc{\good}{\mathrm{good}}
\nc{\tA}{\widetilde\A}
\nc{\Vz}{{V}\kern-1.1ex\raisebox{1.5ex}[0ex][0ex]{$\cdot$}}
\nc{\bxes}[1]{\raisebox{.9ex}{$\cdot$}{\kern#1}\raisebox{0ex}{$\cdot$}
{\kern#1}\raisebox{-.9ex}{$\cdot$}}
\nc{\ssum}{\mathop{\mbox{\normalsize{${\sum}$}}}\limits}
\nc{\ct}{{\mbox{\tiny$\mathrm{CT}$}}}
\nc{\pr}{\mathrm{pr}}
\nc{\qr}{\mathrm{rp}}

\DeclareMathOperator{\ID}{\mathsf{id}}

\DeclareMathOperator{\AUT}{\mathsf{Aut}}

\DeclareMathOperator{\RES}{\mathsf{res}}

\DeclareMathOperator{\PI}{\Pi}

\nc{\Fs}{\ensuremath{\rm F}}

\nc{\isotf}{\overset{
{\rule{0pt}{.9ex}%
{\raisebox{-.6ex}[0ex][-.7ex]{$\mspace{3mu}\sim\mspace{3mu}$}}}}
{\longleftrightarrow}}
\nc{\tN}{\tilde\N}
\nc{\tens}{\mathop\otimes\limits}
\nc{\super}{\mathrm{super}}
\nc{\Mods}{\on{Mod}_\super}
\nc{\rev}{\mathrm{rev}}
\nc{\clif}{\mathfrak{C}}
\nc{\clifm}{\clif^{-}}
\nc{\Fct}{\mathrm{Fct}}
\nc{\Fcts}{\mathrm{Fct}_\super}

\nc{\Ks}{\on{K^\super}}
\nc{\ts}{\widetilde{s}}
\nc{\KUGIRI}{\circ}
\nc{\Sym}{\mathfrak{S}}
\nc{\FF}{\mathcal{F}}
\nc{\BF}{\mathcal{B}}
\nc{\BFF}{\widetilde{\mathcal{B}}}

\nc{\cc}{\mathfrak{c}}
\nc{\SK}{\mathcal{KS}}
\nc{\noi}{\noindent}
\nc{\odd}{{\mathrm{odd}}}
\nc{\even}{{\mathrm{even}}}
\nc{\bs}{\ol{s}}
\nc{\Khc}[1][n]{\ol{\mathcal{KHC}}_{#1}}
\nc{\Ohc}[1][n]{\ol{\mathcal{OHC}}_{#1}}
\nc{\KHC}[1][n]{\mathcal{K}{\mathcal{HC}}_{#1}}
\nc{\OHC}[1][n]{\mathcal{O}{\mathcal{HC}}_{#1}}
\nc{\IODD}{I_{\odd}}
\nc{\IEVEN}{I_{\even}}

\newcommand{\MOD}[1]{{\mathsf{Mod}({#1})}}

\nc{\IRED}{I}

\nc{\MH}{\mathcal{H}}

\newcommand{\MAPSTO}{\longmapsto}

\newcommand{\RP}[1]{\mathsf{RP}_3}

\nc{\bwr}{\mbox{\large$\wr$}}
\nc{\At}[1][{{i,j}}]{\mathscr{A}_{{#1}}}
\nc{\tAt}[1][{{i,j}}]{\widetilde{\mathscr{A}}_{{#1}}}
\nc{\Bt}[1][{{i,j}}]{\mathscr{B}_{{#1}}}
\nc{\tBt}[1][{{i,j}}]{\widetilde{\mathscr{B}}_{{#1}}}
\nc{\prt}[1]{\mathrm{par}(#1)}
\nc{\er}{\mathrm{e}}
\nc{\ec}{\mathrm{e}^-}

\nc{\GCM}{GCM}
\nc{\HCO}{\widehat{\CO}}
\nc{\tCO}{\widetilde{\CO}}
\nc{\red}[1]{{\color{red}{#1}}}
\nc{\T}{\mathbb{T}}
\nc{\HC}{\mathsf{HC}}
\nc{\EV}{\mathsf{ev}}
\nc{\HRC}[1][n]{\widehat{\mathrm{RC}}_{{#1}}}
\nc{\hc}[1][n]{\ol{\mathcal{HC}}_{{#1}}}
\nc{\bphi}{\bar{\phi}}
\nc{\hgt}{\mathrm{ht}}

\begin{document}
%\title{pattern avoidances and multiplicities of maximal weights of affine Lie algebra modules via RSK correspondences}
\title[pattern avoidance in weight multiplicities]{pattern avoidance seen in multiplicities of maximal weights of affine Lie algebra representations}

\author{Shunsuke Tsuchioka}
\address{Graduate School of Mathematical Sciences, University of Tokyo,
Komaba, Meguro, Tokyo, 153-8914, Japan}
\thanks{The first author was supported in part by JSPS Kakenhi Grants 26800005.}
\email{tshun@kurims.kyoto-u.ac.jp}

\author{Masaki Watanabe}
\address{Graduate School of Mathematical Sciences, University of Tokyo,
Komaba, Meguro, Tokyo, 153-8914, Japan}
\email{mwata@ms.u-tokyo.ac.jp}
%\thanks{The author was supported by the EPSRC Postdoctoral Fellowship EP/G050244/1.}

\date{Nov 10, 2015}
\keywords{weight multiplicities, affine Lie algebras, pattern avoidance, maximal weights, Kashiwara's crystal,
Littelmann's path model, RSK correspondence, plane partitions, orbit Lie algebras, quantum binomial coefficients, categorification,
Hecke algebras, symmetric groups, modular representation theory, Mullineux involution}
\subjclass[2010]{Primary~17B67, Secondary~05A05}

\begin{abstract}
We prove that the multiplicities of
certain maximal weights of $\mathfrak{g}(A^{(1)}_{n})$-modules are
counted by pattern avoidance on words.
This proves and generalizes a conjecture of Misra-Rebecca. %~\cite[Conjecture 3.9, Conjecture 4.13]{MR}.
We also prove similar phenomena in types $A^{(2)}_{2n}$ and $D^{(2)}_{n+1}$.
Both proofs are applications of Kashiwara's crystal theory.

%This proves and further generalizes a conjecture of Misra-Rebecca.
%Using a result of Naito-Sagaki on Kashiwara crystals fixed by diagram automorphism, 
%we see similar phenomena in type $A^{(2)}_{2n}$ and $D^{(2)}_{n+1}$.
%occurence of pattern avoidance in type $A^{(2)}_{2n}$ and $D^{(2)}_{n+1}$.
\end{abstract}

\maketitle
%\tableofcontents
\section{Introduction}
Let $\GEE=\GEE(A)$ be a Kac-Moody Lie algebra associated with a symmetrizable GCM $A$.
For each dominant integral weight $\Lambda\in\MP{A}$, we have
the integrable highest weight module $V(\Lambda)$ and
the set of weights $\MYWEIGHT{\Lambda}{A}:=\{\mu\in\HA^*\mid V(\Lambda)_{\mu}\ne 0\}$ with the Weyl group $W$ action.
Studies of the multiplicities of weight spaces, i.e.,
$\MYMULT{\Lambda}{\mu}{A}:=\dim V(\Lambda)_{\mu}$ for $\mu\in \MYWEIGHT{\Lambda}{A}$,
occupy a central position in combinatorial representation theory.
For example, popular algebro-combinatorial ingredients
such as Young Tableaux, Kashiwara's crystal etc
are directly related to such dimension countings.

On the other hand, sometimes
information on $\MYWEIGHT{\Lambda}{A}$ or $\MYMULT{\Lambda}{\mu}{A}$
gives that on representation theory of
seemingly different algebras (and vice-versa) via categorification.
For example,
by virtue of Lascoux-Leclerc-Thibon-Ariki theory and its subsequent developments, 
%when $\GEE=\HSL{p}=\mathfrak{g}(A^{(1)}_{p-1})$ for $p\geq 2$, 
we know
that $\MYWEIGHT{\Lambda}{A^{(1)}_{p-1}}$ parameterizes the blocks of certain
cyclotomic Hecke algebras (a.k.a. Ariki-Koike algebras) $\mathcal{H}$ 
and under this identification it is known that
\bna
\item the orbit space $\MYWEIGHT{\Lambda}{A^{(1)}_{p-1}}/W$ enumerates the possible derived equivalence classes of blocks of $\mathcal{H}$~\cite[\S7.2]{CR},
\item $\MYMULT{\Lambda}{\mu}{A^{(1)}_{p-1}}$ tells us the number of irreducible modules of the block~\cite[Theorem A]{LM}.
\ee
Similar theorems are expected for other types of ``Hecke algebras'', such as KLR algebras, Hecke-Clifford algebras etc
by choosing $A$ suitably.

A rough structure of $\MYWEIGHT{\Lambda}{A}$ is known when $A$ is affine.

\Prop[{\cite[\S12.6]{Kac}}]\label{kaccite}
Let $A$ be an affine GCM.
For $\Lambda\in\MP{A}$, we have
\begin{align*}
\MYWEIGHT{\Lambda}{A}=\bigsqcup_{\lambda\in\MYMAX{\Lambda}{A}}\{\lambda-n\delta\mid n\geq 0\}
\end{align*}
where $\MYMAX{\Lambda}{A}$ is the set of all maximal weights of $V(\Lambda)$ defined as follows.
\begin{align*}
\MYMAX{\Lambda}{A}=\{\lambda\in \MYWEIGHT{\Lambda}{A}\mid \lambda+\delta\not\in \MYWEIGHT{\Lambda}{A}\}.
\end{align*}
\label{Fact1}
\enprop

Clearly, $\MYMAX{\Lambda}{A}$ is $W$-invariant (i.e., $\MYMAX{\Lambda}{A}=W\cdot(\MYMAX{\Lambda}{A}\cap\MP{A})$).
It is known that the set of dominant maximal weights $\MYMAX{\Lambda}{A}\cap\MP{A}$ is finite~\cite[Proposition 12.6]{Kac}.

When $\Lambda$ is level 1, the Hecke algebras appearing in the the aforementioned correspondence via categorification 
%when $\Lambda$ is level 1
%to the case %$\GEE=A^{(1)}_{p-1}$ and 
%of level 1 $\Lambda$ 
are Iwahori-Hecke algebras of type A. Note that
$\MYMAX{\Lambda}{X}\cap\MP{X}=\{\Lambda\}$ when $\Lambda$ is level 1 and $X$ is affine A,D,E type~\cite[Lemma 12.6]{Kac}.
In a course of a study of representation theory of Iwahori-Hecke algebras of type B,
the first author studied the set of dominant maximal weights $\MYMAX{\Lambda_0+\Lambda_s}{A^{(1)}_{p-1}}\cap\MP{A^{(1)}_{p-1}}$ for $0\leq s<p$.
%level 2 weight structures of of $\HSL{p}=\mathfrak{g}(A^{(1)}_{p-1})$.
%The maximal weights and their multiplicities are given as follows.

\Def
Let $p\geq 2$ be an integer (not necessarily a prime).
For $\ell\geq 1$ and $t,u$ with 
$t\geq 0, \ell+t<p-\ell+1$ and 
$u\leq p, \ell<u-\ell+1$, we define 
2 elements of the root lattice $Q$ of $\HSL{p}=\GEE(A^{(1)}_{p-1})$ as follows.
\begin{align*}
\lambda^p_{\ell,t} =  \ell\alpha_0 +
\begin{pmatrix}
\ell\alpha_1 + \cdots \ell\alpha_t \\
+(\ell-1)\alpha_{t+1} + (\ell-2)\alpha_{t+2} + \cdots + \alpha_{\ell+t-1} \\
+\alpha_{p-\ell+1} + \cdots + (\ell-2)\alpha_{p-2} + (\ell-1)\alpha_{p-1}
\end{pmatrix}, \\
\mu^p_{\ell,u} = \ell\alpha_0 +
\begin{pmatrix}
(\ell-1)\alpha_{1} + (\ell-2)\alpha_{2} + \cdots + \alpha_{\ell-1} \\
+\alpha_{u-\ell+1} + \cdots + (\ell-2)\alpha_{u-2} + (\ell-1)\alpha_{u-1} \\
+\ell\alpha_{u} + \cdots \ell\alpha_{p-1}
\end{pmatrix}
\end{align*}
\edf

Recall that $A=A^{(1)}_{p-1}=(2\delta_{ij}-\delta_{i+1,j}-\delta_{i-1,j})_{i,j\in \Z/p\Z}$ and $I=\Z/p\Z$ (see Figure \ref{zu}).
Throughout, we sometimes identify the set $I=\Z/p\Z$ with $\{0,1,\cdots,p-1\}$.

We note that for $p\geq 2$ and when $t=0,u=p$, 
$\lambda^p_{\ell,0}$ is defined exactly when $\mu^p_{\ell,p}$ is defined and in this case we have $\lambda^p_{\ell,0}=\mu^p_{\ell,p}$.
For a Dynkin diagram automorpshim (see \S~\ref{NaitoSagaki}) $\diagautom:I\isoto I,i\MAPSTO -i$, 
we have $\diagautom(\lambda^p_{\ell,t})=\mu^p_{\ell,p-t},\diagautom(\mu^p_{\ell,u})=\lambda^p_{\ell,p-u}$.

The dominant maximal weights and their multiplicities are given as follows.

\Th[{\cite[Theorem 1.4]{Ts1}}]
Let $p\geq 2$ be an integer (not necessarily a prime) and consider a level 2 weight $\Lambda=\Lambda_0 + \Lambda_s$
of $\HSL{p}$ for some $0\leq s<p$. We have
%The set of all dominant maximal weights $\max(\Lambda)\cap\MP^+$ and their multiplicities are described as follows.
\bna
\item {$\displaystyle \MYMAX{\Lambda}{A^{(1)}_{p-1}}\cap\MP{A^{(1)}_{p-1}}=\{\Lambda\}\sqcup\{\Lambda-\lambda^p_{\ell,s}
\mid 1\leq \ell\leq \lfloor\frac{p-s}{2}\rfloor\}
\sqcup\{\Lambda-\mu^p_{\ell,s}\mid 1\leq \ell\leq \lfloor \frac{s}{2}\rfloor\}$}.
\item $\MYMULT{\Lambda}{\Lambda-\lambda^p_{\ell,s}}{A^{(1)}_{p-1}} = \DDD_{\ell,s}$, $\MYMULT{\Lambda}{\Lambda-\mu^{p}_{\ell,s}}{A^{(1)}_{p-1}} = \DDD_{\ell,p-s}$.
\ee
\enth

Here $\DDD_{n,m}$ is the number of lattice paths
from $(0,0)$ to $(n+m,n)$ with steps $(1,0)$ and $(0,1)$ that does not exceed the diagonal $y=x$.
It is not difficult to see $\DDD_{n,m}=\frac{m+1}{n+m+1}{2n+m\choose n}$~\cite[Exercise 6.20.b]{St2}.
\begin{center}
\includegraphics{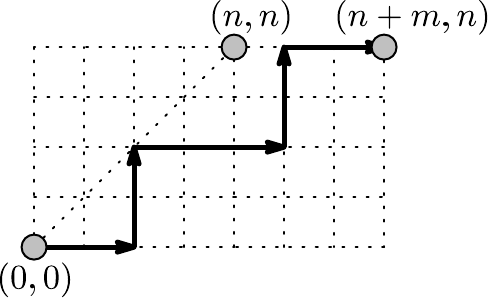}
\end{center}

\begin{figure}
\begin{align*}
\begin{array}{r@{\quad}lrlr@{\quad}l}%rl}
A_1^{(1)}=D_2^{(2)} & \node{}{\alpha_0} \Leftrightarrow \node{}{\alpha_1} &
A_{\ell}^{(1)}  & 
\begin{array}{c}
\raisebox{-12pt}{\rotatebox{30}{$-\!\!-\!\!-$}}\node{}{\alpha_0}\raisebox{0pt}{\rotatebox{-30}{$-\!\!-\!\!-$}} \\[-7pt]
\node{}{\alpha_1}-\node{}{\alpha_2}-\cdots-\node{}{\alpha_{\ell}} 
\end{array} &
D_{\ell+1}^{(2)}  &\node{}{\alpha_0}\Leftarrow \node{}{\alpha_1}-\cdots-\node{}{\alpha_{\ell-1}}\Rightarrow\node{}{\alpha_\ell} 
\\
A_2^{(2)} &\node{}{\alpha_0} {\quad\!\!\!\!\!}\Llleftarrow{\!} \node{}{\alpha_1} &
A_{2\ell}^{(2)}  &\node{}{\alpha_0}\Leftarrow \node{}{\alpha_1}-\cdots-\node{}{\alpha_{\ell-1}}\Leftarrow\node{}{\alpha_{\ell}} &
 {} & {}
\end{array}
\end{align*}
\caption{Affine Dynkin diagrams of A,D,E.}
\label{zu}
\end{figure}

For a higher level $\Lambda\in\MP{A^{(1)}_{p-1}}$, the structure of the set $\MYMAX{\Lambda}{A^{(1)}_{p-1}}\cap\MP{A^{(1)}_{p-1}}$ gets complicated,
but one can easily see the following whose proof will be recalled in \S\ref{inclsubse}.
\Lemma\label{incl}
For $\Lambda=k\Lambda_0+\Lambda_s$ where $k\geq 1$ and $0\leq s<p$, we have
\begin{align*}
\{\Lambda-\lambda^p_{\ell,s}\mid 1\leq \ell\leq \lfloor\frac{p-s}{2}\rfloor\}
\bigsqcup
\{\Lambda-\mu_{\ell,s}\mid 1\leq \ell\leq \lfloor \frac{s}{2}\rfloor\}
\subseteq
\MYMAX{\Lambda}{A^{(1)}_{p-1}}\cap\MP{A^{(1)}_{p-1}}.
\end{align*}
\enlemma

Based on an observation that $\DDD_{n,0}$ is the Catalan number and thus 
the number of $321$-avoiding permutations of $n$~\cite[Exercise 6.19.ee]{St2},
Misra-Rebecca conjectured %the following 
a link between multiplicities of certain maximal weights of $\HSL{p}$-modules and
pattern avoidance.

%On the other hand, $\DDD_{n,0}$ is the number of $321$-avoiding permutations of $n$~\cite[6.19.ee]{St2}.

\Con[{\cite[Conjecture 4.13]{MR}}]\label{mrcon}
For $1\leq \ell\leq\lfloor p/2\rfloor$,
$\MYMULT{(k+1)\Lambda_0}{(k+1)\Lambda_0-\lambda^p_{\ell,0}}{A^{(1)}_{p-1}}$ is equinumerous to 
%the number of 
$((k+2),(k+1),\cdots,2,1)$-avoiding permutations of $\ell$.
\encon

Our main theorem proves and generalizes it in the following way.

\Th\label{ThmA}
Let $p\geq 2$ be an integer and consider a level $k+1$ weight of the form $\Lambda=k\Lambda_0 + \Lambda_s$
of $\HSL{p}$ for some $0\leq s<p$ and $k\geq 1$. 
Then, for $1\leq \ell\leq \lfloor\frac{p-s}{2}\rfloor$, $\MYMULT{\Lambda}{\Lambda-\lambda^p_{\ell,s}}{A^{(1)}_{p-1}}$ is equinumerous to
%the number of 
shuffles of $0^s,1,2,\cdots,\ell$ (there are $s$ zeroes)
that has no strictly decreasing subsequence of length $k+2$.
%sequence $w=(w_1,\cdots,w_{\ell+s})$ such that
%\bnum
%\item for any $0\leq j\leq \ell$, we have $\#\{1\leq i\leq \ell+s\mid w_i=j\}=1+(s-1)\cdot\delta_{0,j}$,
%\item there is no $1\leq i_1<i_2<\cdots<i_{k+2}\leq \ell+s$ such that $w_{i_1}>\cdots>w_{i_{k+2}}$.
%\ee
\enth

By symmetry, for $0<s<p$ and $1\leq\ell\leq \lfloor \frac{s}{2}\rfloor$, $\MYMULT{\Lambda}{\Lambda-\mu^p_{\ell,s}}{A^{(1)}_{p-1}}$ 
is equal to %(when $s=0$, there is no $\ell$ that $\mu^p_{\ell,s}$ makes sense)
\begin{align*}
%\MYMULT{\Lambda}{\Lambda-\mu^p_{\ell,s}}{A^{(1)}_{p-1}}=
\#\{\textrm{shuffles of $0^{p-s},1,2,\cdots,\ell$ that has no strictly decreasing subsequence of length $k+2$}\}.
\end{align*}

Our proof is based on a result of Ariki-Kreiman-Tsuchioka
which characterize the connected component 
(known as Kleshchev multipartitions in 
modular representation theory of Hecke algebras) of $A^{(1)}_{p-1}$-crystal
$B(a\Lambda_0+b\Lambda_s)\subseteq B(\Lambda_0)^{\otimes a}\otimes B(\Lambda_s)^{\otimes b}$ in the tensor 
product~\cite[Corollary 9.6]{AKT}.
This result is a combinatorial incarnation of Littelmann's result~\cite[Theorem 10.1]{Li2}.

While a link between multiplicities of maximal weights of $\HSL{p}$-modules and
pattern avoidance was first observed in ~\cite{MR},
we see similar appearances of pattern avoidance in multiplicities of maximal weights of
affine Lie algebras of types $A^{(2)}_{2n}$ and $D^{(2)}_{n+1}$ (for a crystal-theoretic distinction of
types $A^{(1)}_{n},A^{(2)}_{2n},D^{(2)}_{n+1}$, see ~\cite[\S1]{Ts2}).
In the following, Lie theoretic objects associated with $\AO$ are written with $\check{}$ attached.

\Th\label{ThmB}
Let $p\geq 2$ be an integer and consider a level $k+1$ weight of the form $\LAMBDA=(k+1)\LAMBDAO$
of $\AO=A^{(2)}_{p-1}$ (resp. $D^{(2)}_{1+p/2}$) depending on $p$ being odd (resp. even) where $k\geq 1$ (see Figure \ref{zu}). 
For $1\leq \ell\leq\lfloor p/2\rfloor$, 
\bna
\item\label{bone} $\gamma_{\ell}:=\LAMBDA-\ell \ALPHAO{0}-(\ell-1)\ALPHAO{1}-\cdots-\ALPHAO{\ell-1}\in \MYMAX{\LAMBDA}{\AO}\cap\MP{\AO}$,
\item\label{btwo} $\MYMULT{\LAMBDA}{\gamma_{\ell}}{\AO}$ is equinumerous to %the number of 
$((k+2),(k+1),k,\cdots,1)$-avoiding involutions of $\ell$.
\ee
\enth

Our proof is based on a result of Naito-Sagaki~\cite[Theorem 4.4]{NS1} on Kashiwara's crystals fixed by a diagram automorphism
which is also an application of Littelmann's path model.

When preparing the paper, ~\cite{MR2} appeared on the arXiv which gives a proof of
Conjecture \ref{mrcon}, i.e., the case of $s=0$ of Theorem \ref{ThmA}.
Note that in ~\cite{MR}, Misra-Rebecca also give a conjectural formula~\cite[Conjecture 3.9]{MR}
on the cardinality $\#(\MYMAX{k\Lambda_0}{A^{(1)}_{p-1}}\cap\MP{A^{(1)}_{p-1}})$ to which we will give
a proof in \S\ref{app} using the $q$-Lucas theorem dates back to Gauss.

\vskip 3mm

\noindent{\bf Notations and conventions.} 
We assume that readers are familiar with Kac-Moody Lie algebras and Kashiwara's crystal 
theory (\cite{Kac} and ~\cite{Kas} are standard references).

For integers $a\geq 0$ and $b\geq 1$, we denote by $a\% b$ the remainder of $a$ by $b$, namely
the unique integer $0\leq c<b$ such that $a-c\in b\Z$.
%$\N=\Z_{\geq 0}$ (resp. $\NNN=\Z_{\geq 1}$) means the set of non-negative (resp. positive) integers.

The set of partitions is denoted by $\PAR$ and
the symbol $\emptypartition$ is reserved for the empty partition.
For a partition $\lambda=(\lambda_1,\lambda_2,\cdots)\in\PAR$, % we define %$m_k(\lambda)=|\{i\geq 1\mid \lambda_i=k\}|$ for $k\geq 1$.
we define $|\lambda|=\sum_{i\geq 1}\lambda_i$ and $\ell(\lambda)=\#\{i\geq 1\mid \lambda_i\ne 0\}(=(\TRANS{\lambda})_1)$.
For $n\geq 0$, we put $\PAR(n)=\{\lambda\in\PAR\mid|\lambda|=n\}$.
For $a\geq 0,b\geq 1$, $(a^b)$ is an abbreviation for a partition $c$ such that $c_1=\cdots=c_b=a$.

The symbol $\RPAR{p}$ (resp. $\CRPAR{p}$) stands for the set of $p$-restricted (resp. $p$-core) partitions for $p\geq 2$.
Recall that $\lambda\in\PAR$ is $p$-restricted (resp. $p$-core)
if $\lambda_i-\lambda_{i+1}<p$ for $i\geq 1$ (reps. if there is no removable $p$-hook). 
Note that $\CRPAR{p}\subseteq \RPAR{p}$.

A semistandard tableaux (SST, for short) is a filling of Young diagram by integers
which are weakly increasing along rows and strictly
increasing along columns. 
A column-strict plane partition (CSPP, for short) is a filling of Young diagram by positive integers
which are weakly decreasing along rows and strictly
decreasing along columns. 
For a SST or a CSPP $T$, we denote by $\SHAPE(T)$ the underlying Young diagram.
The content $\CONT(T)$ of $T$ is a multiset of the numbers filled in $T$.

For $\lambda\in\PAR$, we denote by $\SST(\lambda)$ (resp. $\RSST(\lambda)$) the set of SST (resp. CSPP)
of shape $\lambda$.
As usual, $\ST(\lambda)$ (resp. $\RST(\lambda)$) means the set of standard tableaux $T$ (resp. reverse standard tableaux), i.e.,
SST (resp. CSPP) such that $\CONT(T)=\{1,2,\cdots, |\lambda|\}$.%whose contents are exactly $1,2,\cdots,|\lambda|$.

Finally, %for a finite-dimensional algebra $A$ over a field $\corr$,
$\MOD{A}$ means the abelian category of finite-dimensional left $A$-modules and $A$-homomorphisms between them
for a finite-dimensional algebra $A$ over a field $\corr$.
We denote by $\IRR(\MOD{A})$ the set of isomorphism classes of simple objects in $\MOD{A}$.

\vskip 3mm

\noindent{\bf Acknowledgements.} S.T. thanks Satoshi Naito for discussions.%s on results of him jointly with Daisuke Sagaki.

\section{Proof of Theorem \ref{ThmA}}\label{pfA}
%\subsection{Crystal theoretic preliminalies}% from ~\cite{AKT}}
%We apply crystal theory to prove Theorem \ref{ThmA}.
%Let $V$ be a set of permutations of $\underbrace{0,\cdots,0}_{s},1,2,\cdots,\ell$ (there are $s$ zeroes)
%that has no strictly decreasing subsequence of length $k+2$.
In this section, $p,k,\ell,s$ are as in Theorem \ref{ThmA}, i.e., $p\geq 2,k\geq 1,0\leq s<p,1\leq \ell\leq \lfloor(p-s)/2\rfloor$.
We will show that $\MYMULT{\Lambda}{\Lambda-\lambda^p_{\ell,s}}{A^{(1)}_{p-1}}=\#V$ where $\Lambda=k\Lambda_0+\Lambda_s$ and
\begin{align*}
V=\{
\textrm{shuffles of $0^s,1,2,\cdots,\ell$ 
that has no strictly decreasing subsequence of length $k+2$}
\}.
\end{align*}

\subsection{Robinson-Schensted-Knuth correspondence}
Recall the Robinson-Schensted-Knuth correspondence (RSK correspondence, for short) (see ~\cite[\S4]{Ful}).
Fix a multiset $J=\{w_1,\cdots,w_{m}\}\subseteq \Z$ with cardinality $m$ (counted with multiplicity).
RSK correspondence gives a bijection between the set of shuffles (or words) of $w_1,\cdots,w_{m}$
and
\begin{align*}
\bigsqcup_{\lambda\in\PAR(m)}\{(P,Q)\in\SST(\lambda)\times\ST(\lambda)\mid\CONT(P)=J\}.
\end{align*}

We mean by $(P,Q)=\RSK(w)$ that a shuffle (or a word) $w$ maps to a pair of tableaux $(P,Q)$ of the same shape via RSK correspondence.
How RSK correspondence respects ordered subsequences of a shuffle is well-known (see ~\cite[\S3]{Ful}).

\Lemma\label{lenRSK}
Let $(P,Q)=\RSK(w)$ with $\lambda=\SHAPE(P)=\SHAPE(Q)$. Then,
$\ell(\lambda)$ %(resp. $\ell(\lambda)$) 
is the length of the largest strictly decreasing %(resp. weakly increasing) 
subsequence of $w$.
\enlemma

In summary, RSK correspondence gives a bijection between $V$ and $V_1$ where
\begin{align*}
V_1= 
\bigsqcup_{\substack{\lambda\in\PAR(\ell+s) \\ \ell(\lambda)\leq k+1}}
\{(P,Q)\in\SST(\lambda)\times\ST(\lambda)
\mid \CONT(P)=\{0^s,1,\cdots,\ell\}\}.
\end{align*}

Thus, we know that there is a bijection $V\isoto V_2$ where
\begin{align*}
V_2 = 
\bigsqcup_{\substack{\lambda\in\PAR(\ell+s) \\ \ell(\lambda)\leq k+1}}
\{(P,Q)\in\RST(\lambda)\times\RST(\lambda)\mid
\textrm{all $\ell+s,\cdots,\ell+1$ appear in the first row of $P$} 
\}.
\end{align*}

For a permutation $w\in\mathfrak{S}_{m}=\AUT(\{1,\cdots,m\})$, we define $\WORD(w)=w_1\cdots w_{m}$ 
which is a shuffle of $\{1,2,\cdots,m\}$
by $w_i=w(i)$ for $1\leq i\leq m$.
We will use the following well-known symmetry in \S\ref{pB} (see ~\cite[\S4]{Ful}).
\Lemma\label{symRSK}
For $w\in\mathfrak{S}_{m}$ with $\RSK(\WORD(w))=(P,Q)$, we have $\RSK(\WORD(w^{-1}))=(Q,P)$.
\enlemma

\subsection{Kleshchev multipartitions}% from ~\cite{AKT}}
Crystal theoretically, the number $\MYMULT{\Lambda}{\Lambda-\lambda^p_{\ell,s}}{A^{(1)}_{p-1}}$ is
translated as the following counting.
\begin{align*}
%\MYMULT{\Lambda}{\Lambda-\lambda^p_{\ell,s}}{A^{(1)}_{p-1}}=
\#\{b:=\mu\otimes\lambda^{(1)}\otimes\cdots\otimes\lambda^{(k)}\in B(\Lambda_s)\otimes B(\Lambda_0)^{\otimes k}
\mid b\in B(\Lambda),
\WT(b)=\Lambda-\lambda^p_{\ell,s}\}.
\end{align*}
Here $B(\Lambda)$ means the naturally embedded one in $B(\Lambda_s)\otimes B(\Lambda_0)^{\otimes k}$.% (Kleshchev multipartitions).

We adapt Misra-Miwa realization~\cite{MM} for $A^{(1)}_{p-1}$-crystal $B(\Lambda_s)$ for $0\leq s<p$.
We need not know the details of this realization such as the definition of Kashiwara operators.
All we need to know is the following basic things and a result~\cite[Corollary 9.6]{AKT}.
\bnA
\item The underlying set of $B(\Lambda_s)$ is $\RPAR{p}$,
\item For each $\lambda\in B(\Lambda_s)$ and each box $x=(i,j)\in\lambda$ (this means $x$ is the box inside
$\lambda$ located at $i$-th row and $j$-th column), $x$ has the quantity 
$\RES(x)=(s-i+j)+p\Z\in\Z/p\Z$, called
the residue of $x$,
\item\label{defwt} For each $\lambda\in B(\Lambda_s)$, 
\begin{align}
\WT(x)=\WT_s(x):=\Lambda_s-\sum_{i\in\Z/p\Z}\#\{x\in\lambda\mid\RES(x)=i\}\cdot\alpha_i.
\label{wtform}
\end{align}
\ee

\Th[{\cite[Corollary 9.6]{AKT}}]\label{aktcite}
Let $b:=\mu\otimes\lambda^{(1)}\otimes\cdots\otimes\lambda^{(k)}\in B(\Lambda_s)\otimes B(\Lambda_0)^{\otimes k}$.
%in the $A^{(1)}_{p-1}$-crystal. 
Then $b\in B(k\Lambda_0+\Lambda_s)$ (i.e., $b$ is a Kleshchev multipartition)
if and only if $\tau_{(p-s)\%p}(\BASE(\mu))\supseteq\ROOF(\lambda^{(1)})$ and $\BASE(\lambda^{(i)})\supseteq \ROOF(\lambda^{(i+1)})$ for all $1\leq i<k$.
\label{akt2}
\enth

Here $\BASE,\tau_m$~\cite{AKT} where $0\leq m<p$ and $\ROOF$~\cite{KLMW} are explicit maps
\begin{align*}
\begin{cases}
\BASE,\ROOF:\RPAR{p}\longrightarrow\CRPAR{p} \\
\tau_m:\CRPAR{p}\longrightarrow\CRPAR{p}
\end{cases}
\end{align*}
and $\lambda'\supseteq\mu'$ means that $\lambda'$ contains $\mu'$ as Young diagrams.
We need not know the precise definitions of maps $\BASE,\ROOF$ and $\tau_m$, however we need the following minimum.
\bna
\item\label{iden} For a $p$-core partition $\lambda$, we have $\lambda=\BASE(\lambda)=\ROOF(\lambda)$~\cite[Definition 2.5,2.8]{AKT}.
\item For a $p$-core partition $\lambda=(\lambda_1,\cdots,\lambda_a)$, 
we have 
$\tau_{m}(\lambda)=(\nu_1,\cdots,\nu_{a+m})$~\cite[Proposition 9.4]{AKT} where
\begin{align*}
\nu_i=
\begin{cases}
\lambda_i+((p-m)\%p) & (1\leq i\leq m) \\
\min\{\lambda_{i}+((p-m)\%p),\lambda_{i-m}\} & (m<i\leq a) \\
\min\{(p-m)\%p,\lambda_{i-m}\} & (a<i\leq a+m).
\end{cases}
\end{align*}
\ee

Note that $\tau_{0}=\ID_{\CRPAR{p}}$ and
$\tau_{m}(\lambda)=\SHIFT_{(p-m)\%p}(\lambda)\cap (\infty^{m},\lambda)$ where $\SHIFT_{t}(\lambda)=(\lambda_i+t)_{i\geq 1}$ for $\lambda\in\PAR$ and $t\geq 0$.
Of course, $\SHIFT_{t}(\lambda)$ and $(\infty^{m},\lambda)$ are not Young diagrams in the usual sense. But in this section,
infinite Young diagram $\nu$ of these forms
only appears as the form $\nu\supseteq \mu$ for a usual finite Young diagram $\mu\in\PAR$.

\Prop\label{xyzz}
As subsets of $\RPARR{p}{k+1}$, we define
\begin{align*}
X &= \{(\mu,\lambda^{(1)},\cdots,\lambda^{(k)})\mid 
\textrm{$(\JOKEN)$ and
$\tau_{(p-s)\%p}(\BASE(\mu))\supseteq\ROOF(\lambda^{(1)}),1\leq\forall i<k,\BASE(\lambda^{(i)})\supseteq \ROOF(\lambda^{(i+1)})$}\},\\
Y &= \{(\mu,\lambda^{(1)},\cdots,\lambda^{(k)})\in(\CRPAR{p})^{k+1}\mid 
\textrm{$(\JOKEN)$ and
$\tau_{(p-s)\%p}(\mu)\supseteq \lambda^{(1)}\supseteq\cdots\supseteq\lambda^{(k)}$}\},\\
Z &= \{(\mu,\lambda^{(1)},\cdots,\lambda^{(k)})\mid 
\textrm{$(\JOKEN)$ and
$\SHIFT_{s}(\mu)\supseteq \lambda^{(1)}\supseteq\cdots\supseteq\lambda^{(k)}$}\},\\
Z' &= \{(\mu,\lambda^{(1)},\cdots,\lambda^{(k)})\in(\CRPAR{p})^{k+1}\mid 
\textrm{$(\JOKEN)$ and
$\SHIFT_{s}(\mu)\supseteq \lambda^{(1)}\supseteq\cdots\supseteq\lambda^{(k)}$}\},
\end{align*}
where $(\JOKEN)$ means the condition $\WT_{s}(\mu)+\sum_{i=1}^{k}\WT_0(\lambda^{(i)})=\Lambda-\lambda^p_{\ell,s}$. Then, we have $X=Y=Z=Z'$.
\enprop

\Proof
First, observe that $(\ast)$ implies $\mu\subseteq(\ell^{\ell+s})$ and $\lambda^{(i)}\subseteq((\ell+s)^{\ell})$ for $1\leq i\leq k$.
Especially, $(\ast)$ implies $\mu,\lambda^{(1)},\cdots,\lambda^{(k)}\in\CRPAR{p}$.
By (\ref{iden}) above, $X=Y$ and $Z=Z'$.

When $s=0$, it is clear that $Y=Z'$. Assume $0<s<p$.
Note that $\tau_{p-s}(\mu)\supseteq \lambda^{(1)}$ if and only if $\SHIFT_{s}(\mu)\supseteq \lambda^{(1)}$ and
$(\infty^{p-s},\mu)\supseteq \lambda^{(1)}$. By $\lambda^{(1)}\subseteq((\ell+s)^{\ell})$, the latter condition is automatically satisfied.
Thus, we get $Y=Z'$.
\QED

In summary, we now know that $\MYMULT{\Lambda}{\Lambda-\lambda^p_{\ell,s}}{A^{(1)}_{p-1}}=\#Z$.

\Def\label{defbeta}
Let $\{\beta_b\mid b\in\Z\}$ be formal linealy independent elements over $\Z$.
\bna
\item for $p\geq 2$, we define a map (where $\bigoplus_{i\in \Z/p\Z}\Z\alpha_i$ is a root lattice of $\HSL{p}=\mathfrak{g}(A^{(1)}_{p-1})$) by
\begin{align*}
\TOKUMEI_p:\bigoplus_{b\in\Z}\Z\beta_b\longrightarrow \bigoplus_{i\in \Z/p\Z}\Z\alpha_i,\quad
\beta_{b}\longmapsto \alpha_{b+p\Z},
\end{align*}
\item for $\ell\geq 1$ and $s\geq 0$, we define
\begin{align*}
\LAM{p}{\ell}{s}
&=
\beta_{-\ell+1}+\cdots+(\ell-2)\beta_{-2}+(\ell-1)\beta_{-1} \\
&\quad\quad
+\ell\beta_0+\cdots+\ell\beta_{s}+(\ell-1)\beta_{s+1}+(\ell-2)\beta_{s+2}+\cdots+\beta_{\ell+s-1}.
\end{align*}
\ee
\edf

\Cor\label{xyzz2}
We have $Z=Z''$ where 
as subsets of $\RPARR{p}{k+1}$ we define
\begin{align*}
Z'' = \{(\mu,\lambda^{(1)},\cdots,\lambda^{(k)})\in(\CRPAR{p})^{k+1}\mid 
\textrm{$(\JOKEN\JOKEN)$ and $\SHIFT_{s}(\mu)\supseteq \lambda^{(1)}\supseteq\cdots\supseteq\lambda^{(k)}$}
\}
\end{align*}
where $(\JOKEN\JOKEN)$ means the condition $\sum_{(i,j)\in\mu}\beta_{s-i+j}+\sum_{a=1}^{k}\sum_{(i,j)\in\lambda^{(a)}}\beta_{-i+j}=\LAM{p}{\ell}{s}$.
\encor

\Proof
The conditions $0\leq s<p$ and $1\leq \ell\leq \lfloor\frac{p-s}{2}\rfloor$ imply
$\TOKUMEI_p(\LAM{p}{\ell}{s})=\lambda^{p}_{\ell,s}$. Thus, $Z''\subseteq Z$.
The reverse inclusion follows from the fact that for $(\mu,\lambda^{(1)},\cdots,\lambda^{(k)})\in Z$
we have $\mu\subseteq(\ell^{\ell+s})$ and $\lambda^{(i)}\subseteq((\ell+s)^{\ell})$ for $1\leq i\leq k$ as in the proof of Proposition \ref{xyzz}.
\QED

\subsection{Plane partitions} 
Recall that a 2-dimensional array of non-negative integers $\pi=(\pi_{ij})_{i,j\geq 1}$ is a plane partition 
if $\pi_{ij}\geq \pi_{i+1,j},\pi_{i,j+1}$ for $i,j\geq 1$ and
the support $\{(i,j)\in \Z_{\geq 1}\times \Z_{\geq 1}\mid \pi_{ij}> 0\}$ is a finite set.
We denote by $\PP$ the set of plane partitions.
%Note that a CSPP is nothing but a plane partition for which the nonzero elements strictly decrease along columns.

\Def
For a plane partition $\pi$, we define 
\begin{align*}
\WT(\pi)=\sum_{a\geq 1}\sum_{(i,j)\in\pi_{\ast,a}}\beta_{j-i}
\end{align*}
as an element of $\bigoplus_{b\in\Z}\Z\beta_b$
where $\pi_{\ast,j}=(\pi_{1j},\pi_{2,j},\cdots)\in\PAR$.% and $\{\beta_b\mid b\in\Z\}$ is a formal linealy independent elements.
%(for the definition of $\WT_0$, see \eqref{wtform}).
\edf

Clearly, we have (see \eqref{wtform}) 
\begin{align}
\TOKUMEI_p(\WT(\pi))=\sum_{a\geq 1}(\Lambda_0-\WT_0(\pi_{\ast,a})).
\label{remwt}
\end{align}

Recall a famous bijection (that appears most frequently in proving MacMahon plane partition generating functions (see ~\cite[Corollary 7.20.3]{St2}))
\begin{align}
\PI:\bigsqcup_{\lambda\in\PAR}\RSST(\lambda)\times\RSST(\lambda)\isoto\PP.
\label{cspppi}
\end{align}
%summarized in ~\cite[\S7.20]{St2}.
The correspondence $(P,Q)\mapsto\PI(P,Q)$ is briefly described as follows 
(for a detailed explanation including an example, see ~\cite[\S7.20]{St2}).
\begin{quotation}
Let $p^a,q^a\in\PAR$ be the $a$-th columns of $P$ and $Q$. Then, 
the $a$-th column of $\PI(P,Q)$ is a partition given by the Frobenius notation $\rho(p^a,q^a)$.
\end{quotation}

We get Lemma \ref{wtpp}
%see the following 
because we have
\begin{align*}
%\WT_0(\rho(p^j,q^j)) 
\sum_{(i,j)\in\rho(p^a,q^a)}\beta_{j-i}
= \sum_{i\geq 0}\#p^a_{>i}\cdot\beta_{{i}}+\sum_{i< 0}\#q^a_{>-i}\cdot\beta_{{i}}
\end{align*}
where $\#r_{>b}$ is the number of parts of $r$ that is bigger than $b$ for $r\in\{p^a,q^a\}$ and $b\in\Z$.

%from the construction of $\PI$.

\Lemma\label{wtpp}
Let $P,Q\in \RSST(\lambda)$ for some $\lambda\in\PAR$. For $\pi=\PI(P,Q)$, we have
\begin{align*}
\lambda_1 &= \ell(\TRANS{\lambda}) = \max\{j\geq 1\mid \pi_{*,j} \neq \emptypartition\},\\
\WT(\pi) &= \sum_{i\geq 0}\#P_{>i}\cdot\beta_{{i}}+\sum_{i< 0}\#Q_{>-i}\cdot\beta_{{i}}
%\label{wtformula}
\end{align*}
where $\#R_{>i}$ is the number of boxes of $R$ whose number is bigger than $i$ for $R\in\{P,Q\}$.
\enlemma

Note that in the setting of Lemma \ref{wtpp}, we have
\bna
\item the coefficient of $\beta_0$ in $\WT(\pi)$ is $|\lambda|$,
\item $P,Q\in \RST(\lambda)$ if and only if $\WT(\pi) = \beta_{-|\lambda|+1}+2\beta_{-|\lambda|+2}+\cdots+|\lambda|\beta_0+\cdots+2\beta_{|\lambda|-2}+\beta_{|\lambda|-1}$.
\ee

\Prop\label{bij23}
The bijection $\PI$ (see \eqref{cspppi}) gives a bijection
\begin{align*}
\PI\circ\SWAP\circ(\TRANS{(\cdot)}\times\TRANS{(\cdot)})|_{V_2}:V_2\isoto V_3,\quad
(P,Q)\MAPSTO\PI(\TRANS{Q},\TRANS{P})
\end{align*}
where $\beta=\beta_{-\ell-s+1}+2\beta_{-\ell-s+2}+\cdots+(\ell+s)\beta_0+\cdots+2\beta_{\ell+s-2}+\beta_{\ell+s-1}$ and
\begin{align*}
V_3 = \{\pi\in\PP\mid \pi_{\ast,1}\supseteq (s^{\ell+s}),\WT(\pi)=\beta,\pi_{\ast,k+2}=\emptypartition\}.
\end{align*}
\enprop

\Proof
Take $(P,Q)\in V_2$. Because the first column of $\TRANS{P}$ contains $\ell+s,\cdots,\ell+1$,
we see $\PI(\TRANS{Q},\TRANS{P})_{\ast,1}\supseteq (s^{\ell+s})$ by the construction of $\PI$.
Thus, $\PI(\TRANS{Q},\TRANS{P})\in V_3$ by Lemma \ref{wtpp}.

Conversely, take $\pi\in V_3$. Since $\PI$ is a bijection,
there are unique $\lambda\in\PAR$ and $P,Q\in\RSST(\lambda)$ such that $\pi=\PI(\TRANS{Q},\TRANS{P})$.
By Lemma \ref{wtpp}, $|\lambda|=\ell+s$, $\ell(\lambda)\leq k+1$ and
$P,Q\in\RST(\lambda)$.
Observe that $\WT(\pi)=\beta$ implies $\pi_{\ast,1}\subseteq ((\ell+s)^{\ell+s})$. Thus, $(s^{\ell+s})\subseteq \pi_{\ast,1}\subseteq((\ell+s)^{\ell+s})$.
From this, we easily see that all $\ell+s,\cdots,\ell+1$ must appear in the first row of $P$. In other words, $(P,Q)\in V_2$.
\QED

In \S\ref{pB}, we will use a symmetry that obviously follows from the construction of $\PI$.

\Lemma\label{rev}
For $\lambda\in\PAR$ and $P,Q\in\RSST(\lambda)$, put $\pi=\PI(P,Q),\pi'=\PI(Q,P)$.
Then, $\pi'_{\ast,i}=\TRANS{(\pi_{\ast,i})}$ for $i\geq 1$.
\enlemma

\subsection{Proof of Theorem \ref{ThmA}}\label{finalA}
%A proof completes if there is a bijection between $V_3$ and $Z$.
Let us define maps $\ANNO$ and $\ANNOO$ by
\begin{align}
{} &\ANNO:
V_3\longrightarrow Z,\quad
\pi\MAPSTO (\mu,\pi_{\ast,2},\pi_{\ast,3},\cdots,\pi_{\ast,k+1}),\label{situ1}\\
{} &\ANNOO:
Z\longrightarrow V_3,\quad
(\mu',\lambda'^{(1)},\cdots,\lambda'^{(k)})\MAPSTO \pi'
\label{situ2}
\end{align}
where (note that $(s^{\ell+s})\subseteq \pi_{\ast,1}\subseteq((\ell+s)^{\ell+s})$ as in the proof of Proposition \ref{bij23} and
$\mu'\subseteq(\ell^{\ell+s}),\lambda'^{(a)}\subseteq((\ell+s)^{\ell})$  for $1\leq a\leq k$ as in the proof of Corollary \ref{xyzz2})
\bna
\item $\mu=(\nu_1-s,\nu_2-s,\cdots,\nu_{\ell+s}-s)$ for $\nu=\pi_{\ast,1}$,
\item $\pi'_{\ast,a+1}=\lambda'^{(a)}$ for $1\leq a\leq k$ and $\pi'_{\ast,1}=(\mu_1+s,\cdots,\mu_{\ell+s}+s)$.
\ee

In \S\ref{welldef}, 
we show that both $\ANNO$ and $\ANNOO$ are well-defined. This completes the proof
because by construction $\ANNO$ and $\ANNOO$ are mutually inverse each other.

\subsection{Well-definedness of maps $\ANNO$ and $\ANNOO$}\label{welldef}
As a preparation, a direct calculation shows 
\begin{align}
\beta-\beta_{\square}=\LAM{p}{\ell}{s}
\label{wtcalc1}
\end{align}
where $\beta=\sum_{(i,j) \in (\ell+s)^{\ell+s}}\beta_{j-i}$ and $\beta_{\square}=\sum_{(i,j)\in(s^{\ell+s})}\beta_{j-i}$
for $\ell\geq 1,s\geq 0$ (see Definition \ref{defbeta} and Proposition \ref{bij23}).

To prove the well-definedness of $\ANNO$ (resp. $\ANNOO$), it is enough to show
\begin{align*}
\textrm{$\WT_s(\mu)+\sum_{a=1}^{k}\WT_0(\pi_{\ast,a+1})=\Lambda-\lambda^{p}_{\ell,s}$}\quad\textrm{ (resp. $\WT(\pi')=\beta$)}
\end{align*}
in the situation of \eqref{situ1} (resp. \eqref{situ2}). A check for it is shown in \S\ref{wellann1} (resp. \S\ref{wellann2}).

\subsubsection{}\label{wellann1}
By $\Lambda_0-\WT_0(\nu)=(\Lambda_s-\WT_s(\mu))+\sum_{(i,j)\in(s^{\ell+s})}\alpha_{(j-i)+p\Z}$, \eqref{remwt} and $\TOKUMEI_p(\LAM{p}{\ell}{s})=\lambda^{p}_{\ell,s}$,
\begin{align*}
\WT_s(\mu)+\sum_{a=1}^{k}\WT_0(\pi_{\ast,a+1})
=
\Lambda-\TOKUMEI_p(\beta)+\sum_{(i,j)\in(s^{\ell+s})}\alpha_{(j-i)+p\Z} 
=
\Lambda-\TOKUMEI_p(\beta-\beta_{\square})=\Lambda-\lambda^{p}_{\ell,s}.
\end{align*}

\subsubsection{}\label{wellann2}
By Corollary \ref{xyzz2} and \eqref{wtcalc1}, 
\begin{align*}
\WT(\pi')=\beta_{\square}+\sum_{(i,j)\in\mu'}\beta_{(s+j)-i}+\sum_{a=1}^{k}\sum_{(i,j)\in\lambda'^{(a)}}\beta_{j-i}=\beta.
\end{align*}

\section{Proof of Theorem \ref{ThmB}}\label{pB}
%\subsection{Crystal theoretic preliminalies from ~\cite{NS}}
In this section, $p,k,\ell$ are as in Theorem \ref{ThmB}, i.e., $p\geq 2,k\geq 1,1\leq \ell\leq \lfloor p/2\rfloor$.
As in \S\ref{pfA}, we keep identifying $\RPAR{p}$ with $B(\Lambda_0)$ as $A^{(1)}_{p-1}$-crystal through Misra-Miwa realization
and use results in \S\ref{pfA} substituting $s=0$.

\subsection{Mullineux involution}
\Def[{see \cite[6.42]{Mat}}]
For each $b\in B(\Lambda_0)=\RPAR{p}$ of the form $b=\KF{i_j}\cdots\KF{i_1}\emptypartition$ for some $i_1,\cdots,i_j\in\Z/p\Z$, 
$\MUMAP(b)=\KF{-i_j}\cdots\KF{-i_1}\emptypartition$ is well-defined.
\edf

As in ~\cite[Proposition 5.12]{AKT}, there is a crystal morphism $S_h:B(\Lambda_0)\to B(h\Lambda_0)$ for $h\geq 1$ with
certain properties. Let us briefly recall what will be needed.
Under the canonical embedding $B(h\Lambda_0)\hookrightarrow B(\Lambda_0)^{\otimes h}$, we can write $S_h(\lambda)$ of the form
\begin{align}
S_h(\lambda)=\lambda^{(1)}\otimes\cdots\otimes\lambda^{(h)}.
\label{temp2}
\end{align}
Denoting \eqref{temp2} as
\begin{align*}
S_h(\lambda)^{1/h}=(\lambda^{(1)})^{\otimes 1/h}\otimes\cdots\otimes(\lambda^{(h)})^{\otimes 1/h}
\end{align*}
and replacing an occurrence of $(\mu^{\otimes 1/h})^{\otimes k}$ with $\mu^{\otimes k/h}$, we can write 
\begin{align}
S_h(\lambda)^{1/h}=\nu_1^{\otimes a_1}\otimes \nu_2^{\otimes a_2-a_1}\otimes \cdots \otimes \nu_{s}^{\otimes 1-a_{s-1}}.
\label{inflmap}
\end{align}
Here $0<a_1<\cdots<a_{s-1}<1$ in $\Q$ and $\nu_1,\cdots,\nu_s\in\RPAR{p}$ are pairwise distinct.

As in ~\cite[Theorem 5.13]{AKT}, for any $\lambda\in\RPAR{p}$, the right hand side of ~\eqref{inflmap}
%\begin{align*}
%S_h(\lambda)^{1/h}=\nu_1^{\otimes a_1}\otimes\nu_2^{\otimes (a_2-a_1)}\otimes\cdots\otimes\nu_s^{\otimes (1-a_{s-1})}
%\end{align*}
is stable for any sufficiently divisible $h\geq 1$. Furthermore,
\bna
\item $\nu_1,\nu_2,\cdots,\nu_{s}\in\CRPAR{p}$~\cite[Theorem 5.13.(1)]{AKT},
\item $\nu_1\supsetneq\nu_2\supsetneq\cdots\supsetneq\nu_{s}$~\cite[Theorem 5.14]{AKT},
\item $\nu_1=\ROOF(\lambda),\nu_s=\BASE(\lambda)$~\cite[Definition 5.17, Corollary 6.4, Corollary 8.5]{AKT}
\item for any sufficiently divisible $h$, we have (see ~\cite[Proof of Proposition 5.21]{AKT})
\begin{align}
S_h(\MUMAP(\lambda))^{1/h}=(\TRANS{\nu_1})^{\otimes a_1}\otimes(\TRANS{\nu_2})^{\otimes (a_2-a_1)}\otimes\cdots\otimes(\TRANS{\nu_s})^{\otimes (1-a_{s-1})}.
\label{multrans}
\end{align}
\ee

\Cor[{\cite[Proposition 5.21]{AKT}}]\label{contmm}
For any $\lambda\in\RPAR{p}$, we have $\BASE(\MUMAP(\lambda))=\TRANS{\BASE(\lambda)}$ and $\ROOF(\MUMAP(\lambda))=\TRANS{\ROOF(\lambda)}$.
\encor

\Cor\label{contm}
For any $\lambda\in\CRPAR{p}$, we have $\MUMAP(\lambda)=\TRANS{\lambda}$.
\encor

\Rem
The involution $\MUMAP:\RPAR{p}\isoto\RPAR{p}$ is known as Mullineux involution 
in modular representation theory of symmetric groups and Hecke algebras (see ~\cite[\S7]{LLT}).
Under the identification via cellular algebra structure (see ~\cite[3.43]{Mat})
\begin{align*}
\RPAR{p}\isoto\bigsqcup_{n\geq 0}\IRR(\MOD{\mathbb{F}_p\mathfrak{S}_{n}}),\quad
\lambda\MAPSTO \SPECHT{p}{\lambda},
\end{align*}
%due to James~\cite{Jam}, 
%which we don't explain here,
%using cellular algebra structure (see ~\cite{Mat})
Ford-Kleshchev showed $\SPECHT{p}{\lambda}\otimes\SIGN{\mathbb{F}_p}\cong \SPECHT{p}{\MUMAP(\lambda)}$ 
that had known as Mullineux conjecture.
Here $\SIGN{\corr}$ is the sign representation for a field $\corr$.
It is a classical result that $S^{\lambda}_{\mathbb{Q}}\otimes\SIGN{\mathbb{Q}}\cong S^{\TRANS{\lambda}}_{\mathbb{Q}}$ where 
$\{S^{\lambda}_{\mathbb{Q}}\mid \lambda\in\PAR(n)\}=\IRR(\MOD{\Q\mathfrak{S}_n})$ are the classical Specht modules.
~\eqref{multrans} says that choosing an appropriate model of ``Young diagram'', $\textrm{-}\otimes\SIGN{\corr}$ is 
always given by transposition of Young diagram
even over positive characteristics (for Hecke algebras, see ~\cite[\S5]{AKT}).
\enrem

\subsection{Diagram automorphisms and orbit Lie algebras}\label{NaitoSagaki}
Let $A=(a_{ij})_{i,j\in I}$ be a symmetrizable GCM with a corresponding Kac-Moody Lie algebra $\GEE=\GEE(A)$.
A diagram automorphism $\diagautom:I\isoto I$ is a bijection such that $a_{\diagautom(i),\diagautom(j)}=a_{ij}$ for $i,j\in I$.
For a symmetrizable GCM with a diagram automorphism, the orbit Lie algebra $\GEEO=\GEE(\AO)$,
which is again a Kac-Moody Lie algebra, is defined as follows (see ~\cite[\S2.2]{FSS}).
\bnum
\item put $c_{ij}=\sum_{k=0}^{N_j-1}a_{i,\diagautom^k(j)}$ for $i,j\in I$ where $N_i=\#\{\diagautom^k(i)\mid k\in\Z\}$,
\item set $\IO=\{i\in I/\diagautom\mid c_{ii}>0\}$ and $\AO=(\check{a}_{ij}:=2c_{ij}/c_{jj})_{i,j\in\IO}$.
\ee

In our case of $A=A^{(1)}_{p-1}=(2\delta_{ij}-\delta_{i+1,j}-\delta_{i-1,j})_{i,j\in \Z/p\Z}$ and $I=\Z/p\Z$,
we adapt 
\begin{align*}
\diagautom:\Z/p\Z\isoto\Z/p\Z,\quad
i\MAPSTO -i
\end{align*}
as a diagram automorphism. Then, the orbit Lie algebra is $\GEEO=\GEE(\AO)$ 
where $\AO=A^{(2)}_{p-1}$ (resp. $D^{(2)}_{1+p/2}$) depending on $p$ being odd (resp. even). Recall that
\begin{align}
\DEO=\begin{cases}
2\ALPHAO{0}+\cdots+2\ALPHAO{(p-3)/2}+\ALPHAO{(p-1)/2} & (\AO=A^{(2)}_{p-1}) \\
\ALPHAO{0}+\cdots+\ALPHAO{p/2} & (\AO=D^{(2)}_{1+p/2}).
\end{cases}
\label{delo}
\end{align}

We identify the set $\{i\in I/\diagautom\mid c_{ii}>0\}$ above with
\begin{align*}
\IO=\begin{cases}
\{0,1,\cdots,(p-1)/2\} & (\AO=A^{(2)}_{p-1}) \\
\{0,1,\cdots,p/2\} & (\AO=D^{(2)}_{1+p/2}).
\end{cases}
\end{align*}

For $i\in\IO$, a direct calculation shows 
\begin{align*}
c_{ii}=\begin{cases}
1 & (\textrm{$p\equiv 1\!\!\!\!\pmod{2}$ and $i=(p-1)/2$}) \\
2 & (\textrm{otherwise}).
\end{cases}
\end{align*}

As in Theorem \ref{ThmB}, Lie theoretic objects associated with $\GEEO$ are written with $\check{}$ attached.

\subsection{Naito-Sagaki's fixed points crystals}
%For $n\geq 1$,
Let $B_n$ be the connected component in $\RPARR{p}{n}\cong B(\Lambda_0)^{\otimes n}$
that is isomorphic to $B(n\Lambda_0)$ as $\GEE$-crystal for $n\geq 1$.
By Theorem \ref{aktcite},
\begin{align*}
B_n=\{
(\lambda^{(1)},\cdots,\lambda^{(n)})\in\RPARR{p}{n}\mid
1\leq\forall i<n, \BASE(\lambda^{(i)})\supseteq \ROOF(\lambda^{(i+1)})
\}.
\end{align*}

By virtue of Naito-Sagaki~\cite[Theorem 4.4]{NS1}, 
the set of fixed points $B_n^{\MUMAP^{n}}$ has a $\GEEO$-crystal structure that is isomorphic to $B(n\LAMBDAO)$.
All we need is the correspndence on weights: 
\begin{quotation}
the weight $\WTO(b)$ of $b=(x_1,\cdots,x_{n})\in B_n^{\MUMAP^n}$ as a $\GEEO$-crystal
is given by 
\begin{align}
\WTO(b)=n\LAMBDAO-\sum_{i\in \IO}m_i\ALPHAO{i}\Longleftrightarrow
\sum_{i=1}^{n}\WT_0(x_i)=n\Lambda_0-\sum_{i\in \IO}\frac{2m_i}{c_{ii}}\sum_{r=1}^{N_{i}-1}\alpha_{\iota(\diagautom^r(i))}
\label{wto}
\end{align}
where $\iota:\IO\hookrightarrow I,i\mapsto i+p\Z$ is an injection (see also ~\cite[(1.2.2)]{NS2}).
\end{quotation}

Since $\ell-1<(p-1)/2$ (resp. $\ell-1<p/2$) for odd $p$ (resp. even $p$),
the right hand side of \eqref{wto} is equal to $(k+1)\Lambda_0-\lambda^{p}_{\ell,0}$
whenever the left hand side of \eqref{wto} is given by 
\begin{align*}
\gamma_{\ell}=(k+1)\LAMBDAO-\ell \ALPHAO{0}-(\ell-1)\ALPHAO{1}-\cdots-\ALPHAO{\ell-1}
\end{align*}
for $n=k+1$. Thus, we have 
$\MYMULT{(k+1)\LAMBDAO}{\gamma_{\ell}}{\AO}=\#(Z'^{\MUMAP^{k+1}})$ where (see Proposition \ref{xyzz})
\begin{align*}
Z' = \{(\lambda^{(1)}\supseteq\cdots\supseteq\lambda^{(k+1)})\in(\CRPAR{p})^{k+1}\mid 
\sum_{i=1}^{k+1}\WT_0(\lambda^{(i)})=(k+1)\Lambda_0-\lambda^p_{\ell,0} \}.
\end{align*}

\subsection{Proof of Theorem \ref{ThmB}}
In \S\ref{finalA}, we presented bijections
\begin{align*}
V_2\isoto
V_3\isoto Z=Z',\quad
(P,Q)\MAPSTO
\pi:=\PI(\TRANS{Q},\TRANS{P})\MAPSTO (\pi_{\ast,1},\cdots,\pi_{\ast,k+1})
\end{align*}
where $V_2=\bigsqcup_{\substack{\lambda\in\PAR(\ell) \\ \ell(\lambda)\leq k+1}}\RST(\lambda)^2$ and
$V_3 = \{\pi\in\PP\mid \WT(\pi)=\sum_{(i,j) \in \ell^\ell} \beta_{j-i},\pi_{\ast,k+2}=\emptypartition\}$.

By Corollary \ref{contm} and Lemma \ref{rev}, we have 
$\#(Z'^{\MUMAP^{k+1}})=\sum_{{\lambda\in\PAR(\ell) \\ \ell(\lambda)\leq k+1}}\#\RST(\lambda)$.
This is equal to $\sum_{\lambda\in\PAR(\ell),\ell(\lambda)\leq k+1}\#\ST(\lambda)$ and
it is equinumerous to $((k+2),(k+1),k,\cdots,1)$-avoiding involution of $\ell$ 
by Lemma \ref{lenRSK} and Lemma \ref{symRSK}.
This completes the proof of Theorem \ref{ThmB} (\ref{btwo}).

We now know that $\MYMULT{(k+1)\LAMBDAO}{\gamma_{\ell}}{\AO}>0$. Thus, to prove Theorem \ref{ThmB} (\ref{bone}), 
it is enough to show that $\MYMULT{(k+1)\LAMBDAO}{\gamma_{\ell}+\DEO}{\AO}=0$ (see Proposition \ref{kaccite}). This follows from Proposition \ref{Fact3} and
the condition $\ell-1<(p-1)/2$ (resp. $\ell-1<p/2$) when $p$ is odd (resp. $p$ is even).

\Prop[{\cite[Proposition 12.5.(a)]{Kac}}]
Let $A$ be an affine GCM.
For $\Lambda\in\MP{A}$, 
\begin{align*}
\MYWEIGHT{\Lambda}{A} =
W\cdot\{\lambda\in\MP{A}\mid\lambda\leq\Lambda\}.
\end{align*}
\label{Fact3}
\enprop

\section{Appendix: on the number of $\MYMAX{k\Lambda_0}{A^{(1)}_{p-1}}\cap\MP{A^{(1)}_{p-1}}$}\label{app}
We prove a conjecture of Misra-Rebecca on the number $\#(\MYMAX{k\Lambda_0}{A^{(1)}_{p-1}}\cap\MP{A^{(1)}_{p-1}})$.
%for $k\geq 1$ and $p\geq 2$.

\Prop[{\cite[Conjecture 3.9]{MR}}]\label{saigo}
for $k\geq 1$ and $p\geq 2$.
\begin{align*}
\#(\MYMAX{k\Lambda_0}{A^{(1)}_{p-1}}\cap\MP{A^{(1)}_{p-1}})
=
\frac{1}{p+k}\sum_{d\Z\supseteq k\Z,p\Z}
\EUPHI(d){(p+k)/d \choose k/d}
\end{align*}
where $\EUPHI(d)=\#(\Z/d\Z)^{\times}$ is Euler's totient function.
\enprop

\subsection{Proof of Lemma \ref{incl}}\label{inclsubse}
Recall that $p\geq 2,k\geq 1,0\leq s<p$ and $\Lambda=k\Lambda_0+\Lambda_s$.
Depending on $s\ne 0$ or not, we define the set $S^{(p,s)}_{k}$ as follows.
\begin{align*}
S^{(p,0)}_{k} &= \{(x_i)_{i=0}^{p}\in \Z^{p+1}\mid x_0=x_p=0,x_1+x_{p-1}\leq k,0<\forall i<p,-x_{i-1}+2x_i-x_{i+1}\geq 0\}, \\
S^{(p,s)}_{k} &= \{(x_i)_{i=0}^{p}\in \Z^{p+1}\mid x_0=x_p=0,x_1+x_{p-1}\leq k-1,0<\forall i<p,\delta_{i,s}-x_{i-1}+2x_i-x_{i+1}\geq 0\}.
\end{align*}

As in ~\cite[\S3.1,\S3.2]{Ts1}, the following gives a bijection.
\begin{align*}
S^{(p,s)}_{k}\isoto \MYMAX{\Lambda}{A^{(1)}_{p-1}}\cap\MP{A^{(1)}_{p-1}},\quad
(x_0,\cdots,x_p)\MAPSTO \Lambda+\sum_{i=0}^{p-1}(x_i+q_0)\alpha_{i}
\end{align*}
where $q_0=\max\{q\leq 0\mid \textrm{$1\leq \forall i<p,x_i+q\leq 0$ and $1\leq \exists i<p,x_i+q=0$}\}$.
Clearly $S^{(p,s)}_{k}\subseteq S^{(p,s)}_{k+1}$ and $q_0$ does not depend on $k$, thus we deduce Lemma \ref{incl}.

\subsection{$q$-binomial coefficients and $q$-Lucas theorem}
Let ${a \brack b}=[a]!/([b]![a-b]!)$ be a $q$-binomial coefficient for $0\leq b\leq a$ and $[c]!=\prod_{n=1}^{c}(q^n-1)/(q-1)$.

\Prop[{\cite[pp.66]{St1}}]\label{qbinomyng}
For any $j,k\geq 0$, we have $\sum_{\substack{\lambda\in\PAR \\ \ell(\lambda)\leq j, \lambda_1\leq k}}q^{|\lambda|}={k+j \brack j}$.
%\begin{align*}
%\sum_{\substack{\lambda\in\PAR \\ \ell(\lambda)\leq j, \lambda_1\leq k}}q^{|\lambda|}={k+j \brack j}.
%\end{align*}
\enprop

%Sagan proved the following congruent property as $q$-Lucas theorem.
The following congruent property for $q$-binomial coefficients is known as $q$-Lucas theorem 
(see also ~\cite[Exercise 14 of Chapter 1]{St1} for Lucas theorem for binomial coefficients).
\Prop[{\cite[Theorem 2.2]{Sag}}]\label{sagthm}
Let $\zeta$ be a primitive $d$-th root of unity where $d\geq 1$. For any $n,j\geq 0$, %we have 
\begin{align*}
{n \brack j}\Big|_{q=\zeta}={\lfloor n/d\rfloor \choose \lfloor j/d\rfloor}\cdot{n\%d \brack j\%d}\Big|_{q=\zeta}.
\end{align*}
\enprop

\subsection{Proof of Proposition \ref{saigo}}
As in \S\ref{inclsubse}, $\#(\MYMAX{k\Lambda_0}{A^{(1)}_{p-1}}\cap\MP{A^{(1)}_{p-1}})=\#S^{(p,0)}_{k}$.

Let us define sets $T$ and $U$ as follows.
\begin{align*}
%S &= \{(x_0,\cdots,x_p)\in \Z^{p+1}\mid x_0=x_p=0,x_1+x_{p-1}\leq k,\forall 0<i<n,-x_{i-1}+2x_i-x_{i+1}\geq 0\},\\
T &= \{(y_1,\cdots,y_p)\in \Z^{p}\mid y_1\geq\cdots\geq y_p,y_1+\cdots+y_p=0,y_1-y_p\leq k\},\\
U &= \{(\lambda_1,\cdots,\lambda_{p-1})\in \Z^{p-1}\mid k\geq\lambda_1\geq\cdots\lambda_{p-1}\geq 0,\lambda_1+\cdots+\lambda_{p-1}\in p\Z\}.
\end{align*}

The following maps are bijections.
\begin{align*}
S^{(p,0)}_{k} &\isoto T,\quad (x_0,\cdots,x_p)\MAPSTO (x_1-x_0,\cdots,x_p-x_{p-1}),\\
T &\isoto U,\quad (y_1,\cdots,y_p)\MAPSTO (y_1-y_p,\cdots,y_{p-1}-y_{p}).
\end{align*}

%By Proposition \ref{qbinomyng}, we have
%\begin{align*}
%\sum_{k\geq\lambda_1\geq\cdots\lambda_{p-1}\geq 0}q^{\lambda_1+\cdots+\lambda_{p-1}}={k+p-1 \brack p-1}
%\end{align*}
%where ${a \brack b}=\frac{[a]!}{[b]![a-b]!}$ is a $q$-binomial coefficient for $0\leq b\leq a$ and $[c]!=\prod_{n=1}^{c}\frac{q^n-1}{q-1}$.

By Proposition \ref{qbinomyng}, we have
\begin{align*}
\#U=\frac{1}{p}\sum_{\zeta^p=1}{k+p-1 \brack p-1}\Big|_{q=\zeta}.
\end{align*}

Let $\zeta$ be a primitive $d$-th root of unity for some $1\leq d\leq p$ with $d\Z\ni p$. Then, %by Proposition \ref{sagthm}, we have
\begin{align}
{k+p-1 \brack p-1}\Big|_{q=\zeta}={\lfloor (k+p-1)/d\rfloor \choose \lfloor (p-1)/d\rfloor}{(k+p-1)\%d \brack d-1}\Big|_{q=\zeta}
\label{temp1}
\end{align}
by Proposition \ref{sagthm}.
The right hand side of ~\eqref{temp1} vanishes unless $k+p-1\equiv d-1\pmod{d}\Leftrightarrow d\Z\ni k$.
When $d\Z\ni k$, the right hand side of ~\eqref{temp1} 
becomes ${\lfloor (k+p-1)/d\rfloor \choose \lfloor (p-1)/d\rfloor}=\binom{(k+p)/d-1}{p/d-1}$.
Thus, we know that $\#(\MYMAX{k\Lambda_0}{A^{(1)}_{p-1}}\cap\MP{A^{(1)}_{p-1}})=\#S^{(p,0)}_{k}=\#U$ is equal to
\begin{align*}
\frac{1}{p}\sum_{d\Z\supseteq k\Z,p\Z}
\EUPHI(d){(p+k)/d-1 \choose p/d-1} 
= \frac{1}{p+k}\sum_{d\Z\supseteq k\Z,p\Z}
\EUPHI(d){(p+k)/d \choose k/d}.
\end{align*}

\end{document}